\documentclass[11pt]{article}
\usepackage{epsfig,makeidx,amsmath,amsfonts,latexsym,subfigure}
\newtheorem{theorem}{Theorem}[section]
\newtheorem{prop}{Proposition}[section]

\newcommand{\E}{{\mathbb E}}

\newcommand {\PP}{\mathbb P}

\begin{document}

\title{Random networks with preferential growth and vertex death}

\author{Maria Deijfen \thanks{Department of Mathematics, Stockholm University, 106 91 Stockholm. Email: mia@math.su.se.}}

\date{July 2010}

\maketitle

\thispagestyle{empty}

\begin{abstract}

\noindent A dynamic model for a random network evolving in continuous time is defined where new vertices are born and existing vertices may die. The fitness of a vertex is defined as the accumulated in-degree of the vertex and a new vertex is connected to an existing vertex with probability proportional to a function $b$ of the fitness of the existing vertex. Furthermore, a vertex dies at a rate given by a function $d$ of its fitness. Using results from the theory of general branching processes, an expression for the asymptotic empirical fitness distribution $\{p_k\}$ is derived and analyzed for a number of specific choices of $b$ and $d$. When $b(i)=i+\alpha$ and $d(i)=\beta$ -- that is, linear preferential attachment for the newborn and random deaths -- then $p_k\sim k^{-(2+\alpha)}$. When $b(i)=i+1$ and $d(i)=\beta(i+1)$, with $\beta<1$, then $p_k\sim (1+\beta)^{-k}$, that is, if also the death rate is proportional to the fitness, then the power law distribution is lost. Furthermore, when $b(i)=i+1$ and $d(i)=\beta(i+1)^\gamma$, with $\beta,\gamma<1$, then $\log p_k\sim -k^\gamma$ -- a stretched exponential distribution. The momentaneous in-degrees are also studied and simulations suggest that their behaviour is qualitatively similar to that of the fitnesses.

\vspace{0.5cm}

\vspace{0.5cm}

\end{abstract}

\section{Introduction}

Empirical studies on real networks have revealed that many of them exhibit features that are not captured by the classical Erd\"{o}s-R\'{e}nyi graph. In particular, many networks tend to have a quite heavy tailed degree distribution, often described by a power law, that is, the fraction of vertices with degree $k$ decays as $k^{-\tau}$ for some exponent $\tau$. To capture this, a number of new graph models have been developed that allow for this type of degree distribution; see for instance the seminal paper by Bollob\'{a}s et al.\ \cite{BJR} (the setup described there includes many previously studied models as special cases) and the references therein.

A model type that has received a lot of attention is based on preferential attachment: new vertices are successively added to the network and are connected to existing vertices proportionally to a function $f(k)=k^\gamma$ of the degree. The mechanism was introduced in the context of network modeling by Barab\'asi and Albert \cite{BA}, who observed that, for $\gamma=1$ it seemed to lead to power law distributions with exponent $\tau=3$ for the degrees. This was proved rigorously in \cite{BRST}. Krapivsky and Redner \cite{KR} later derived heuristic results indicating that the regimes $\gamma\leq 1$ and $\gamma>1$ are qualitatively different: for $\gamma>1$ the degree distribution is degenerate in the sense that there is a single dominant vertex that is connected to almost every other vertex. This was confirmed in \cite{OS}. Recently Dereich and M\"{o}rters \cite{DM} has studied the degree evolution in the sublinear regime $\gamma<1$.

Most preferential attachment models in the literature are defined in discrete time. Rudas et al.\ \cite{RT} however studies the model from a new perspective, embedding it in continuous time and interpreting it as a continuous time branching process. This opens up for the use of well established results from the theory of branching processes. Compared to the traditional techniques for analyzing models with linear preferential attachment, the approach has the advantage that it applies for any $\gamma\leq 1$. We mention also the paper by Athreya et al.\ \cite{AGS}, where the linear case $\gamma=1$ is studied by aid of a different type of branching process embedding.

The purpose of the current paper is to define and analyze a preferential attachment model embedded in continuous time where vertices are not just added but may also die, and where the death rate of a vertex is taken to be a function of its previous success, quantified by its accumulated degree. Using results from the theory of general branching processes, a characterization of the asymptotic distribution of the accumulated degrees is derived and analyzed in various special cases. The momentaneous degrees are also discussed and analyzed by means of simulations. This generalizes results from \cite{RT}. More specifically, the generalization compared to \cite{RT} lies in that we work with a birth and death process instead of a pure birth process. We mention also the paper \cite{BL} by Britton and Lindholm, where a model is treated that is similar in spirit but not based on preferential attachment.

Preferential attachment models incorporating deaths of vertices are much less studied than pure growth models, while real networks may well be affected by vertex emigration. Existing rigorous results are restricted to the case when vertices die randomly, that is, the death mechanism does not depend on features of the vertices (age, degree etc); see \cite{CL} and \cite{CFV}. See also \cite{MGN} and \cite{DC} for heuristic results, the latter concerning a situation where vertices die with a probability inversely proportional to their degree. Concentration results are so far comparably weak, see e.g.\ \cite[Lemma 4.1]{CFV}, and are lacking for models with more complicated death mechanisms.

\subsection{Definition of the model}

We now proceed to define the model. As already described, a vertex population evolving in continuous time is considered where existing vertices give birth to new vertices but may also die. More specifically, a vertex that is alive at time $t$ and that has given birth to $i$ new vertices during its life time, gives birth to an $i+1$:st vertex at rate $b(i)$, that is, the probability that the vertex gives birth in the time interval $(t,t+\varepsilon)$ is $b(i)\varepsilon+o(\varepsilon)$. Furthermore, the vertex dies at rate $d(i)$, that is, its probability of dying during $(t,t+\varepsilon)$ is $d(i)\varepsilon+o(\varepsilon)$. Let $\xi$ be a random point process on $[0,\infty)$ with points (representing the births) generated according to the above rules. Starting from a single root vertex at time 0, each vertex $x$ reproduces according to an independent copy $\xi_x$ of $\xi$. This means that the vertex population evolves as a birth and death process in continuous time and hence fits into the framework of general branching processes; see e.g.\ Jagers \cite[Chapter 6]{J}. Below we shall formulate conditions that guarantee that the process is supercritical, so that the vertex population indeed grows to infinity with positive probability.

Write $\mathcal{Z}_t$ for the set of vertices that are alive at time $t$, and, for $x\in\mathcal{Z}_t$, let $A_t^x$ denote the total number of children that vertex $x$ has produced up to time $t$ (including children that are no longer alive), also referred to as the \emph{fitness} of the vertex. The time until the next birth in the population is exponentially distributed with parameter $\sum_{x\in\mathcal{Z}_t} b\left(A_t^x\right)$ and the time until the next death is exponentially distributed with parameter $\sum_{x\in\mathcal{Z}_t} d\left(A_t^x\right)$. Hence the probability that the next event that occurs in the population is a birth is given by
\begin{equation}\label{eq:bprop}
\frac{\sum_{x\in\mathcal{Z}_t} b\left(A_t^x\right)}{\sum_{x\in\mathcal{Z}_t} (b\left(A_t^x\right)+ d\left(A_t^x\right))},
\end{equation}
and the probability that the birth occurs at vertex $x$ is proportional to $b(A_t^x)$. We shall see later that, under certain assumption on $b$ and $d$, the ratio (\ref{eq:bprop}) converges almost surely to a number in $(0,1)$.

The vertex population described above naturally gives rise to a network structure: directed edges are created between a vertex and its children, pointing from the child to the mother. When a vertex dies it turns into a \emph{ghost} and all edges associated with it are turned into \emph{ghost links}. According to the dynamics, the rate at which a vertex produces/attracts new vertices depends on its previous success, measured by its in-degree, including ghost links. Weather or not it is realistic to include the ghost links depends on the setting. In many situations it is indeed natural to assume that the ability of a vertex to acquire new edges is affected by the total historical success of the vertex rather than just the contributions generated by vertices that are still present in the network (think e.g.\ of coauthor-ship networks, sexual networks, recruitment networks etc). A model where ghost links are excluded in the definition of the fitness can be analyzed by the same means as the present model. The calculations however become cumbersome and do not lead to explicit expressions. We elaborate on this and present simulation results in Section 5.

\subsection{Results}

Let $\rho(t)$ be the density of $\xi$, that is,
$$
\rho(t)=\lim_{\varepsilon\to 0}\varepsilon^{-1}\PP((t,t+\varepsilon)\mbox{ contains a point from }\xi)
$$
and define the Laplace transform
\begin{equation}\label{eq:rhohat}
\hat{\rho}(\lambda)=\int_0^\infty e^{-\lambda t}\rho(t)dt.
\end{equation}
Clearly $\hat{\rho}(\lambda)$ increases as $\lambda\to 0$, and we write
\begin{equation}\label{eq:rhoinf}
\underline{\lambda}=\inf\{\lambda>0:\hat{\rho}(\lambda)<\infty\}.
\end{equation}
Our first result is an expression for $\hat{\rho}(\lambda)$.

\begin{prop}\label{prop:rhohat}
We have that
$$
\hat{\rho}(\lambda)=\sum_{k=0}^\infty\prod_{i=0}^k\frac{b(i)}{\lambda+b(i)+d(i)}.
$$
\end{prop}

\noindent Throughout, we shall assume that the functions $b$ and $d$ are such that
\begin{equation}\tag{C.1}
\hat{\rho}(\lambda^*)=1\quad\mbox{has a solution }\lambda^*\in(0,\infty).
\end{equation}
In fact, we shall assume that
\begin{equation}\tag{C.2}
\lambda^*>\underline{\lambda}.
\end{equation}
It is well-known from the theory of branching processes that, when $\underline{\lambda}<\infty$, the condition (C.1) is equivalent to the process being supercritical. To be more precise, write $\xi(t)$ for the $\xi$-measure of $[0,t]$ and let $Z_t=|\mathcal{Z}_t|$ be the number of vertices alive at time $t$. Then $\PP(Z_t\to\infty)>0$ if and only if $\E[\xi(\infty)]>1$ if and only if $\lambda^*>0$; see \cite[Chapter 6]{J}. Writing
\begin{equation}\label{eq:elog}
\hat{\rho}(\lambda)=\sum_{k=0}^\infty e^{\sum_{i=0}^k\log\left(1-\frac{\lambda+d(i)}{\lambda+b(i)+d(i)}\right)}
\end{equation}
and using the Taylor expansion $\log(1-x)=-x+O(x^2)$ reveals that a necessary condition for $\hat{\rho}(\lambda)<\infty$ is that
$$
\sum_{i=0}^\infty\frac{\lambda+d(i)}{\lambda+b(i)+d(i)}=\infty,
$$
and hence a necessary condition for $\underline{\lambda}<\infty$ (which is indeed necessary for (C.1)) is that
$$
\sum_{i=0}^\infty\frac{1+d(i)}{b(i)+d(i)}=\infty.
$$
Furthermore, a sufficient condition for having $\hat{\rho}(\lambda)<\infty$ for all $\lambda>0$ (so that $\underline{\lambda}=0$) is that
\begin{equation}\label{eq:rhohatfinite}
\frac{1+d(i)}{d(i)+b(i)}\geq ci^{-\gamma} \quad\mbox{for some }\gamma<1.
\end{equation}
Indeed, a bound of the same order then applies for $\frac{\lambda+d(i)}{\lambda+b(i)+d(i)}$ and hence, by Taylor expansion and integral approximation, the exponent in (\ref{eq:elog}) decays at least like $-k^{1-\gamma}$, making the sum convergent. Also note that (C.1) -- or, more specifically, the fact that $\underline{\lambda}<\infty$ -- guarantees that the process does not explode in finite time: a sufficient condition for non-explosion is that $\E[\xi(t)]<\infty$ for some $t>0$, see \cite[Theorem 6.2.2]{J}, and, since $\rho(t)dt=\E[\xi(dt)]$, it is clear that $\underline{\lambda}$ cannot be finite if this does not hold.

We are now ready to formulate the main theorem, which is a characterization of the asymptotic fitness distribution in the vertex population. To this end, recall that $Z_t$ is the number of alive vertices at time $t$ and write $Z_t^k$ for the number of alive vertices at time $t$ with fitness $k$.
Furthermore, let $L$ be a random variable distributed as the life time of a vertex.

\begin{theorem}\label{th:dist}
Assume that $b$ and $d$ satisfy conditions (C.1) and (C.2). Then, for all $k$, on the event $\{Z_t\to\infty\}$, almost surely
$$
\lim_{t\to\infty}\frac{Z_t^k}{Z_t}=\frac{C}{\lambda^*+b(k)+d(k)}
\prod_{i=0}^{k-1}\frac{b(i)}{\lambda^*+b(i)+d(i)}
$$
where $C^{-1}=\int_0^\infty e^{-\lambda^*t}\PP(L>t)dt$ and the empty product arising for $k=0$ is defined as 1.
\end{theorem}

The rest of the paper is organized as follows. In Section 2 we state the results we shall need from the theory of branching process and, in Section 3, Proposition \ref{prop:rhohat} and Theorem \ref{th:dist} are proved. Then, in Section 4, the limiting distribution in Theorem \ref{th:dist} is analyzed in more detail for a number of specific choices of $b(i)$ and $d(i)$. In Section 5, a formula is given for the the asymptotic degree distribution in a network where the birth and death rates are based on the momentaneous in-degrees rather than the fitnesses (that is, the accumlated in-degrees). Unfortunately this expression is hard to analyze, but simulation results are presented that indicate that the qualitative behaviour of the distribution is the same as for the fitnesses. Finally, some possibilities for further work are mentioned in Section 6.

\section{Preliminaries on branching processes}

The corpus of literature on branching processes is vast and we do not intend to give a survey; for this, see Jagers \cite{J} or the more expository book by Haccou et al.\ \cite{HJV}. We shall however briefly explain the concepts that are needed to establish Theorem \ref{th:dist}. These revolve around general branching processes counted with random characteristics. In a general branching process, the individuals reproduce independently of each other according to a random point process $\xi$ on $[0,\infty)$ during a life-time $L$, where there are no restrictions on the dependence structure between $\xi$ and $L$. An individual $x$ can be represented as $x=(j_1,\ldots,j_n)$ indicating that $x$ is the $j_n$:th child of the $j_{n-1}$:th child $\ldots$ of the $j_1$:th child of the ancestor. Hence the space of possible individuals is
$$
\Omega=\bigcup_{i=0}^\infty \mathbb{N}^n.
$$
Write $\sigma_x$ for the time when the individual $x$ is born, that is, $\sigma_0=0$ and, if $x$ is the $i$:th child of the individual $x'$, then $\sigma_x=\sigma_{x'}+\inf\{t:\xi_{x'}(t)\geq i\}$.

Roughly, when counting a branching processes with a random characteristic, a random score is assigned to each individual and the quantity of interest is then the total score in the population. More specifically, let $\phi(t)$ be a real-valued stochastic process with $\phi(t)=0$ for $t\leq 0$ giving the score of an individual of age $t$. The score of the individual $x$ at time $t$ is then given by $\phi_x(t)=\phi(t-\sigma_x)$ and the total score in the population at time $t$ becomes
$$
Z_t^\phi=\sum_{x\in \Omega}\phi_x(t).
$$
Examples of common characteristics are $\phi(t)=\mathbf{1}\{t\leq 0\}$, in which case $Z_t^\phi$ is the number of individuals born up to time $t$, and $\phi(t)=\mathbf{1}\{0\leq t< L\}$, in which case $Z_t^\phi$ is the number of individuals that are alive at time $t$ -- this quantity we shall use the notation $Z_t$ for. These characteristics are both individual, meaning that $\xi_x$ depends only on the process of descendants of $x$. Characteristics that involve information also about preceding generations are called relational. An example of such a characteristic that we shall be interested in below is when an individual is assigned the value 1 if its mother is still alive and 0 otherwise.

The Laplace transform $\hat{\rho}(\lambda)$ of the density $\rho(t)$ of $\xi$ and the related quantity $\underline{\lambda}$ are defined as in (\ref{eq:rhohat}) and (\ref{eq:rhoinf}). The solution $\lambda^*$ to the equation $\hat{\rho}(\lambda)=1$ is called the Malthusian parameter and the process is subcritical, critical or supercritical depending on weather $\lambda^*<0$, $=0$ or $>0$. A more frequent formulation of the definition of $\lambda^*$ in the literature is as the solution to the equation
$$
\int_0^\infty e^{-\lambda t}\E[\xi(dt)]=1,
$$
where, as in Section 1.1, $\xi(t)=\xi([0,t])$. Since $\E[\xi(dt)]=\rho(t)dt$, the formulations are clearly equivalent.

Now, Theorem \ref{th:dist} will turn out to be a direct consequence of the following result concerning almost sure convergence of ratios of random characteristics.

\begin{theorem}[Nerman (1981)]\label{th:nerman} Consider a supercritical branching process with Malthusian parameter $\lambda^*>\underline{\lambda}$ and let $\phi$ and $\psi$ be two individual characteristics satisfying
\begin{equation}\label{eq:charcond}
\E[\sup_t\{e^{-\lambda t}\phi(t)\}]<\infty\quad\mbox{for some }\lambda<\lambda^*
\end{equation}
and likewise for $\psi$. Then, on the event $\{Z_t\to\infty\}$, almost surely
$$
\frac{Z_t^\phi}{Z_t^\psi}\to\frac{\hat{\phi}(\lambda^*)}{\hat{\psi}(\lambda^*)}
$$
where $\hat{\phi}(\lambda)=\int_0^\infty e^{-\lambda t}\E[\phi(t)]dt$.
\end{theorem}

\noindent \textbf{Remark 1.} The characteristics we shall work with will in general be bounded and, for such characteristics, (\ref{eq:charcond}) is trivially satisfied.\medskip

\noindent \textbf{Remark 2.} As pointed out in \cite[pp 163]{HJV}, the conclusion of the theorem holds true also for relational characteristics provided that these depend only on a finite number of preceding generations; see \cite{JN} for details.

\section{Proofs}

We now return to the specific setting defined in Section 1.1, where $\xi$ is the point process generated by the described birth and death rules.\medskip

\noindent \textbf{Proof of Proposition \ref{prop:rhohat}.} The expression for $\hat{\rho}(\lambda)$ follows from straightforward calculations. First note that
$$
\hat{\rho}(\lambda)=\sum_{k=0}^\infty b(k)\int_0^\infty e^{-\lambda t}\PP(\xi(t)=k\cap L>t)dt,
$$
where trivially
$$
\PP(\xi(t)=k\cap L>t)=\PP(\xi(t)=k)-\PP(\xi(t)=k\cap L\leq t).
$$
Consider a generic vertex (the root vertex for instance) and write $T_i$ for the time, counting from the birth of the vertex, when its $i$:th child is born, with $T_i:=\infty$ if $\xi(\infty)<i$. Also, when $T_i<\infty$, define $\Delta^b_i=T_{i+1}-T_i$ and let $\Delta^d_i$ be the time, counting from $T_i$, until the vertex dies. For $k=0$, we obtain
\begin{eqnarray*}
\int_0^\infty e^{-\lambda t}\PP(\xi(t)=0\cap L>t)dt & = & \int_0^\infty e^{-(\lambda +b(0)+d(0)) t}dt\\
& = & \frac{1}{\lambda+b(0)+d(0)},
\end{eqnarray*}
where $e^{-(b(0)+d(0))t}$ is the probability that nothing happens in the process up to time $t$. To deal with $k\geq 1$, note that
$$
\PP(\xi(t)=k) = \PP(\xi(t)=k|T_k<\infty)\PP(T_k<\infty),
$$
with
\begin{eqnarray*}
\PP(\xi(t)=k|T_k<\infty) & = & \PP(\xi(t)\geq k|T_k<\infty)\\
& & -\PP(\xi(t)\geq k+1|T_{k+1}<\infty)\PP(T_{k+1}<\infty|T_k<\infty).
\end{eqnarray*}
Here
$$
\PP(\xi(t)\geq k|T_k<\infty)=\PP\left(\sum_{i=0}^{k-1}\Delta^b_i\leq t|T_k<\infty\right),
$$
where, conditionally on $T_k<\infty$, the variables $\Delta^b_0,\ldots,\Delta^b_{k-1}$ are independent and $\Delta^b_i$ is exponentially distributed with parameter $b(i)+d(i)$. Similarly,
\begin{eqnarray*}
\PP(\xi(t)=k\cap L\leq t) & = & \PP(\xi(t)=k\cap L\leq t|T_k<\infty\cap T_{k+1}=\infty)\cdot\\
&& \PP(T_k<\infty\cap T_{k+1}=\infty)
\end{eqnarray*}
with
$$
\PP(\xi(t)=k\cap L\leq t|T_k<\infty\cap T_{k+1}=\infty)=\PP\left(\sum_{i=0}^{k-1}\Delta^b_i+\Delta^d_i\leq t\big|T_k<\infty\cap T_{k+1}=\infty\right),
$$
where, conditionally on $T_k<\infty\cap T_{k+1}=\infty$, the variables $\Delta^b_1,\ldots,\Delta^b_{k-1},\Delta^d_k$ are independent, $\Delta^b_i$ exponentially distributed with  parameter $b(i)+d(i)$ and $\Delta^d_k$ exponentially distributed with parameter $b(k)+d(k)$. Also note that
$$
\PP(T_k<\infty)=\prod_{i=0}^{k-1}\frac{b(i)}{b(i)+d(i)}=:q_k
$$
$$
\PP(T_{k+1}<\infty|T_k<\infty)=\frac{b(k)}{b(k)+d(k)}=:q_k'
$$
$$
\PP(T_k<\infty\cap T_{k+1}=\infty)=\frac{d(k)}{b(k)+d(k)}\prod_{i=0}^{k-1}\frac{b(i)}{b(i)+d(i)}=:q_k''.
$$
For independent exponential variables $\{X_i\}$, where $X_i$ has parameter $b(i)+d(i)$, write
$$
I_k=\int_0^\infty e^{-\lambda t}\PP(X_0+\ldots+X_{k-1}\leq t)dt.
$$
Combining all of the above, we get that
\begin{equation}\label{eq:nastan_klart}
\hat{\rho}(\lambda)=\frac{b(0)}{\lambda+b(0)+d(0)}+\sum_{k=1}^\infty b(k)\Big(q_kI_k-(q_kq_k'+q_k'')I_{k+1}\Big),
\end{equation}
and all that remains is to note that,
\begin{eqnarray*}
I_k & = & \frac{1}{\lambda}\int_0^\infty e^{-\lambda t}f_{X_0+\ldots+X_{k-1}}(t)dt\\
& = & \frac{1}{\lambda}\prod_{i=0}^{k-1}\frac{b(i)+d(i)}{\lambda+b(i)+d(i)},
\end{eqnarray*}
where the first equality follows from partial integration and the second from independence. Substituting this, and the expressions for $q_k$, $q_k'$ and $q_k''$, into (\ref{eq:nastan_klart}) gives the desired formula for $\hat{\rho}(\lambda)$.\hfill$\Box$\medskip

We proceed to prove the main result, Theorem \ref{th:dist}, which, as already mentioned, follows from Theorem \ref{th:nerman}.\medskip

\noindent \textbf{Proof of Theorem \ref{th:dist}.} Consider the characteristic
$$
\phi(t)=\mathbf{1}\{\xi(t)=k\cap 0<t<L\},
$$
that is, $Z_t^\phi=Z_t^k$ is the number of alive vertices with $k$ children at time $t$, and let $\psi(t)=\mathbf{1}\{0<t<L\}$ so that $Z_t^\psi=Z_t$ is the number of alive vertices at time $t$. These characteristics are both bounded, and hence, since $b$ and $d$ are assumed to be such that $(C.1)$ and $(C.2)$ are satisfied, all assumptions of Theorem \ref{th:nerman} are fulfilled. It follows that, on the event $Z_t\to\infty$, the ratio $Z_t^k/Z_t$ converges almost surely to the limit identified in Theorem \ref{th:nerman}. The expression in the nominator becomes
$$
\hat{\phi}(\lambda^*)=\int_0^\infty e^{-\lambda^* t}\PP(\xi(t)=k\cap L>t)dt,
$$
and identical calculations as in the proof of Proposition \ref{prop:rhohat} yields the formula in Theorem \ref{th:dist}.\hfill$\Box$

\section{Examples}

We now investigate the behavior of the limiting fitness distribution in Theorem \ref{th:dist} for a few specific choices of birth and death rates $b$ and $d$. The limiting distribution is henceforth denoted by $\{p_k\}$. Furthermore, for two sequences $\{a_k\}$ and $\{b_k\}$ of real numbers, we write $a_k\sim b_k$ to denote that $a_k/b_k$ converges to a strictly positive constant as $k\to\infty$.\medskip

\noindent \textbf{Example 1.} First consider the case when $b(i)=i+\alpha$ and $d(i)=\beta$ with $\alpha>0$, that is, the birth rate is a linear function of the fitness and the vertices die at a constant rate. This means that a new vertex is connected to an existing vertex with a probability that depends linearly on the fitness of the existing vertex, and the vertices die at random. Using the relation
$$
\prod_{i=0}^k(i+c)=\frac{\Gamma(k+1+c)}{\Gamma(c)},
$$
we get that
$$
\hat{\rho}(\lambda)=\sum_{k=0}^\infty \prod_{i=0}^k\frac{i+\alpha}{\lambda+i+\alpha+\beta}=\frac{\Gamma(\lambda+\alpha+\beta)}{\Gamma(\alpha)}\sum_{k=0}^\infty \frac{\Gamma(k+1+\alpha)}{\Gamma(k+1+\alpha+\beta+\lambda)}.
$$
Van der Hofstad et al. \cite[Lemma 5.4]{HMH} derives the formula
\begin{equation}\label{eq:piet}
\sum_{k=0}^n\frac{\Gamma(k+a)}{\Gamma(k+c)}=\frac{1}{1+a-c}\left(\frac{\Gamma(n+1+a)}{\Gamma(n+c)}-
\frac{\Gamma(a)}{\Gamma(c-1)}\right),
\end{equation}
and, by Stirlings formula, $\Gamma(k+c)/\Gamma(k)\sim k^c$ as $k\to \infty$. It follows that
$$
\hat{\rho}(\lambda) = \left\{ \begin{array}{ll}
                      \infty & \mbox{if $\lambda\leq 1-\beta$};\\
                      \frac{\alpha}{\lambda+\beta-1} & \mbox{otherwise}.
                    \end{array}
            \right.
$$
Hence $\lambda^*=1+\alpha-\beta$, implying that the process is supercritical if $\beta<\alpha+1$. For the fitness distribution we get that
$$
p_k=\frac{C}{k+1+2\alpha}\prod_{i=0}^{k-1}\frac{i+\alpha}{i+1+2\alpha}=C'\frac{\Gamma(k+\alpha)}{\Gamma(k+2+2\alpha)}\sim k^{-(2+\alpha)}.
$$
The limiting behavior of the fitness distribution is hence unaffected by the vertex death. This is explained by the fact that the vertices die randomly, that is, given that the next thing to occur is a vertex death, the vertex to die is chosen uniformly at random, and such a mechanism does not affect the tail behavior of the power law distribution generated by the birth mechanism. With $\alpha=1$ the expression for $p_k$ is in agreement with previous results for pure linear preferential attachment.

Recall that, given the development of the process up to time $t$, the probability that the next thing occurring is a vertex birth is given by (\ref{eq:bprop}), where $A_t^x$ is the fitness of vertex $x$ at time $t$. If the characteristics given by $\phi_x(t)=b(A_{t-\sigma_x}^x)$ and $\psi_x(t)=d(A_{t-\sigma_x}^x)$ satisfy the condition (\ref{eq:charcond}), then it follows from Theorem \ref{th:nerman} that the probability converges almost surely to a number in $(0,1)$. This number is in general difficult to determine, but in the current case, denoting the total fitness $A_t^{tot}=\sum_{x\in\mathcal{Z}_t}A_t^x$, we get
$$
\frac{\sum_{x\in\mathcal{Z}_t} b\left(A_t^x\right)}{\sum_{x\in\mathcal{Z}_t} (b\left(A_t^x\right)+ d\left(A_t^x\right))}=\frac{A_t^{tot}+\alpha Z_t}{A_t^{tot}+(\alpha+\beta)Z_t}.
$$
The total fitness in the alive population can be obtained by counting the population by the characteristic $\psi(t)=\mathbf{1}\{t>0$ and the mother of the individual is alive$\}$. Hence, by Theorem \ref{th:nerman} and the second remark following it, the ratio $A_t^{tot}/Z_t$ converges to 1 almost surely and it follows that
$$
\lim_{t\to\infty}\PP(\mbox{the next event is a vertex birth})=\frac{1+\alpha}{1+\alpha+\beta}.
$$
This probability is strictly larger than 1/2 precisely when $\beta<\alpha+1$, that is, when the process is supercritical.\medskip

\noindent \textbf{Example 2.} Next consider the case when $b(i)=i+1$ and $d(i)=\beta(i+1)$, that is, both the birth and the death rate are linear functions of the fitness. In this case, by (\ref{eq:rhohatfinite}), we have $\hat{\rho}(\lambda)<\infty$ for all $\lambda>0$, and it is given by
\begin{eqnarray*}
\hat{\rho}(\lambda) & = & \sum_{k=0}^\infty \prod_{i=0}^k\frac{i+1}{\lambda+(1+\beta)i+(1+\beta)}\\
& = & \Gamma\left(\frac{\lambda+1+\beta}{1+\beta}\right)\sum_{k=0}^\infty
\frac{\Gamma(k+2)}{\Gamma\left(k+1+\frac{\lambda+1+\beta}{1+\beta}\right)}\frac{1}{(1+\beta)^k}\\
& = & F\left(1,1,\frac{\lambda+1+\beta}{1+\beta}, \frac{1}{1+\beta}\right)-1,
\end{eqnarray*}
where $F$ denotes Gauss hypergeometric function; see Abramowitz and Stegun \cite[Chapter 15]{AS}. The equation $\hat{\rho}(\lambda^*)=1$ cannot be solved explicitly but a numerical solution is easily obtained for a given value of $\beta$. The process is supercritical precisely when $\beta<1$, since the expected number of children produced by a vertex is given by
\begin{equation}\label{eq:exp1}
\E[\xi(\infty)] = \sum_{k=1}^\infty \PP(\xi(\infty)\geq k)= \sum_{k=1}^\infty\prod_{i=0}^{k-1}\frac{b(i)}{b(i)+d(i)}= \sum_{k=1}^\infty \frac{1}{(1+\beta)^k}=\frac{1}{\beta}.
\end{equation}
The asymptotic fitness distribution is
\begin{eqnarray*}
p_k & = & \frac{C}{\lambda^*+(1+\beta)(k+1)}\prod_{i=0}^{k-1}\frac{i+1}{\lambda^*+(1+\beta)(i+1)}\\
& = & C'\frac{\Gamma(k+1)}{\Gamma\left(k+3+\frac{\lambda^*}{1+\beta}\right)}\frac{1}{(1+\beta)^k}\\
& \sim &k^{-\left(2+\frac{\lambda*}{1+\beta}\right)}\frac{1}{(1+\beta)^k}.
\end{eqnarray*}
A mechanism where the vertices die at a rate proportional to their fitness hence destroys the power law distribution obtained without vertex deaths regardless of how small the proportionality constant is and instead produces a (super)exponentially decaying fitness distribution.

We remark that $b(i)=(i+1)^\gamma$ and $d(i)=\beta(i+1)^\gamma$ gives rise to the same exponential decay for any $\gamma>0$. Indeed,
\begin{eqnarray*}
p_k & = &  \frac{1}{\lambda^*+(1+\beta)(k+1)^\gamma}\prod_{i=0}^{k-1}
\frac{(i+1)^\gamma}{\lambda^*+(1+\beta)(i+1)^\gamma}\\
& \sim & \frac{1}{(1+\beta)^k}\,e^{-\sum_{i=0}^{k-1}\log\left(1+\frac{\lambda^*}{(1+\beta)(i+1)^\gamma}\right)},
\end{eqnarray*}
where the sum in the exponent grows like $k^{1-\gamma}$ for $\gamma<1$ and is convergent for $\gamma>1$.\medskip

\noindent \textbf{Example 3.} Next take $b(i)=i+1$ and $d(i)=\beta(i+1)^\eta$ for $\beta<1$ and $\eta\in(0,1)$. Also in this case $\hat{\rho}(\lambda)<\infty$ for all $\lambda$ and, by comparison with the previous example, $\beta<1$ guarantees that the process is supercritical. For the fitness distribution we get that
$$
\log p_k\sim -\sum_{i=0}^k\log\left(1+\frac{\lambda^*+\beta(i+1)^\eta}{i+1}\right)\sim -k^\eta.
$$
Vertex deaths determined by a sublinear function of the fitness hence also destroy the power law distribution obtained without vertex death, bringing it down to a stretched exponential distribution. In contrast to the case when the vertices die at random however, the tail behavior of the distribution is changed by the vertex deaths.

Again we remark that $b(i)=(i+1)^\gamma$ and $d(i)=\beta (i+1)^\eta$ for arbitrary $\gamma>\eta>0$ with $\gamma-\eta\in(0,1)$ gives rise to the same type of fitness distribution: a similar derivation as above implies that $\log p_k\sim -k^{1-(\gamma-\eta)}$.\medskip

\noindent \textbf{Example 4.} Finally, consider a situation with $b(i)=\alpha$ and $d(i)=(i+1)^{-1}$, that is, new vertices attach randomly to existing vertices and vertex deaths occur inversely proportional to the fitness. Analogously to (\ref{eq:exp1}) we get that
$$
\E[\xi(\infty)]=\Gamma\left(1+\frac{1}{\alpha}\right)\sum_{k=1}^\infty
\frac{\Gamma(k+1)}{\Gamma\left(k+\frac{1}{\alpha}\right)},
$$
and using (\ref{eq:piet}) it follows that
$$
\E[\xi(\infty)] = \left\{ \begin{array}{ll}
                      \infty & \mbox{if $\alpha\geq 1$};\\
                      \frac{\alpha}{1-\alpha} & \mbox{otherwise}.
                    \end{array}
            \right.
$$
Hence the process is supercritical if and only if $\alpha>1/2$. The Laplace transform is
\begin{eqnarray*}
\hat{\rho}(\lambda) & = & \sum_{i=0}^\infty \left(\frac{\alpha}{\lambda+\alpha}\right)^{k+1}
\frac{\Gamma(k+2)}{\Gamma\left(k+2+\frac{1}{\lambda+\alpha}\right)}\\
& = & F\left(1,1,1+\frac{1}{\lambda+\alpha},\frac{\alpha}{\lambda+\alpha}\right)-1,
\end{eqnarray*}
where again $F$ is the hypergeometric function, and the equation $\hat{\rho}(\lambda^*)=1$ is easily solved numerically for a given value of $\alpha$. As for the fitness distribution, we have
\begin{eqnarray*}
p_k & = & \frac{1}{\lambda^*+\alpha+\frac{1}{k+1}}\prod_{i=0}^{k-1}\frac{\alpha}{\lambda^*+\alpha+\frac{1}{k+1}}\\
& = & \frac{1}{\alpha}\left(\frac{\alpha}{\lambda^*+\alpha}\right)^{k+1}
\frac{\Gamma(k+2)}{\Gamma\left(k+2+\frac{1}{\lambda^*+\alpha}\right)} \\
& \sim & \left(\frac{\alpha}{\lambda^*+\alpha}\right)^k k^{-\frac{1}{\lambda^*+\alpha}},
\end{eqnarray*}
that is, the fitness distribution decays (super)exponentially at rate $\alpha/(\lambda^*+\alpha)$. The fact that vertices that are not successful in acquiring edges dies at a higher rate is hence not sufficient in itself to explain the power law distributions for the fitnesses.

\section{The in-degrees}

As mentioned in Section 1.1, the fitnesses (that is, the accumulated in-degrees) are not seldom the primary object of interest when studying a network incorporating vertex death. However, there are also situations when the momentaneous in-degrees, henceforth referred to as just the in-degrees, are the relevant quantities. This is presumably the case for instance when modeling the World Wide Web. In this section we demonstrate how a formula for the asymptotic in-degree distribution follows from the same method as for the fitnesses and present simulation results that show that the in-degree distribution behaves in a similar way as the fitness distribution.

Consider a vertex population that evolves in the same way as described above, with the difference that the birth and death rates are based on the in-degrees rather than the accumulated in-degrees. More specifically, a vertex that is alive at time $t$ and has $i$ \emph{alive} children at time $t$ gives birth to a new vertex at rate $b(i)$ and dies at rate $d(i)$. Let $\tilde{\xi}$ be a random point process on $[0,\infty)$ with points (representing the births) generated according to these rules. We use the same notation as before equipped with a wave-hat for quantitites related to $\tilde{\xi}$, omitting the original hat for the Laplace transform of the reproduction function. The Laplace transform is given by
\begin{equation}\label{eq:tilde_rho}
\tilde{\rho}(\lambda)=\sum_{k=0}^\infty b(k)\int_0^\infty\PP\left(\tilde{\xi}(t)=k\cap \tilde{L}>t\right)dt,
\end{equation}
and $\tilde{\lambda}^*$ is defined by $\tilde{\rho}(\tilde{\lambda}^*)=1$. Furthermore, if $b(i)$ and $d(i)$ satisfies the analogues of (C.1) and (C.2), it follows from Theorem \ref{th:nerman} that
\begin{equation}\label{eq:in_deg}
\lim_{t\to\infty}\frac{\tilde{Z}_t^k}{\tilde{Z}_t}=\frac{\int_{0}^\infty e^{t\tilde{\lambda}^*}\PP\left(\tilde{\xi}(t)=k\cap \tilde{L}>t\right)dt}{\int_0^\infty\PP(\tilde{L}>t)dt},
\end{equation}
that is, the fraction of vertices with in-degree $k$ converges almost surely and the limit is given by the right-hand of (\ref{eq:in_deg}). Unfortunately it does not seem possible to derive a more explicit expression for the limit in this case. The event that a vertex of age $t$ has $k$ alive children  occurs if the vertex has given birth to $n$ children for some $n\geq k$ and $n-k$ of these children have died before the vertex has reached the age $t$. To find an expression for the probability of this event that leads to a tractable formula for the integral in (\ref{eq:tilde_rho}) and in the numerator of (\ref{eq:in_deg}) does not seem possible.

To get a grip of the asymptotic in-degree distribution, we turn to simulations. These indicate that the qualitative behavior of the in-degrees is the same as for the fitnesses. Figure 1 shows simulated in-degrees and fitnesses plotted on a log-log scale for $b(i)=i+1$ and $d(i)=0.5$. In Example 1 it was shown that the asymptotic fitnesses in this case follow a power-law distribution with exponent 3, and Figure 1 reveals that the in-degrees seem to follow the same type of distribution, although the exponent appears to be slightly larger. Figure 2 shows simulated in-degrees and fitnesses for $b(i)=i+1$ and $d(i)=0.5(i+1)$, with a log-scale on the $y$-axis. The asymptotic fitness distribution was shown in Example 2 to have an exponentially decaying tail in this case, and the same seems to be true for the in-degrees. Finally, Figure 3 shows simulated in-degrees and fitnesses on a log-log scale for $b(i)=i+1$ and $d(i)=0.5(i+1)^{0.5}$. The plot reveals a slightly curved behavior for both the in-degrees and the fitnesses which is significative for a stretched exponential distribution and in agreement with the analytical result for the fitness distribution. Note that, in all three cases, the tail of the in-degree distribution seems to decay faster than the tail of the fitness distribution. All simulations have been run until the number of alive vertices equals 100,000.

\section{Further work}

There are a number of questions about the studied model that deserve further investigation. What about the global properties of the obtained structure for instance? If the ghost links are included one gets a random tree, analyzed in Rudas and Toth \cite{RT}, and excluding the ghost links gives rise to a forest. Will there be an infinite component in this forest? How does the answer depend on $b(i)$ and $d(i)$? Another issue to look into is the time evolution of the process. For models without vertex death it is known that the maximal degree grows polynomially of the same order as the first vertex; see M\'{o}ri \cite{Mo}. Is this the case also in the current setting for the maximal fitness? Recently Dehrich and M\"{o}rters \cite{DM} have shown for a model with attachment function $k^\gamma$ and no vertex death that there is a phase transition at $\gamma=1/2$: if $\gamma>1/2$ there is a persistent hub that will be of maximal degree at all but at most finitely many times, while, if $\gamma<1/2$, there is no persistent hub. Is there a similar phenomenon in a model with vertex death?\bigskip

\noindent \textbf{Acknowledgement.} I thank Remco van der Hofstad for drawing my attention to continuous time branching processes as a possible tool for analyzing preferential attachment models with vertex death, and Gerard Hooghiemstra for pointing me to Lemma 5.4 in \cite{HMH}.

\begin{figure}[p]
\begin{center}
\mbox{{\epsfig{file=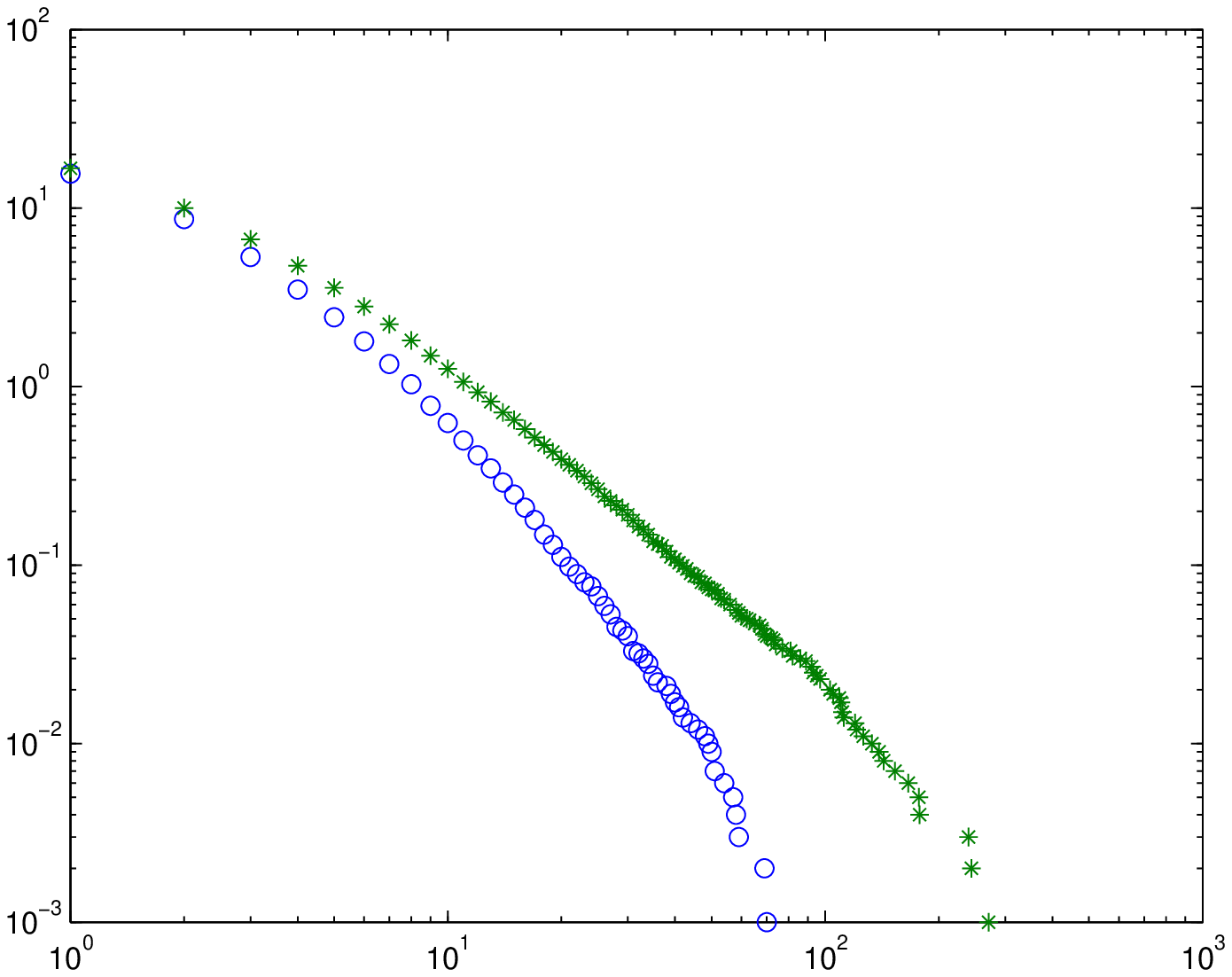,width=0.8\textwidth,height=0.55\textheight}}}
\end{center}
\caption{Simulation with $b(i)=i+1$ and $d(i)=0.5$. The figure shows the tail probabilities of the fitness distribution (stars) and the in-degree distribution (circles).}
\end{figure}

\begin{figure}[p]
\begin{center}
\mbox{{\epsfig{file=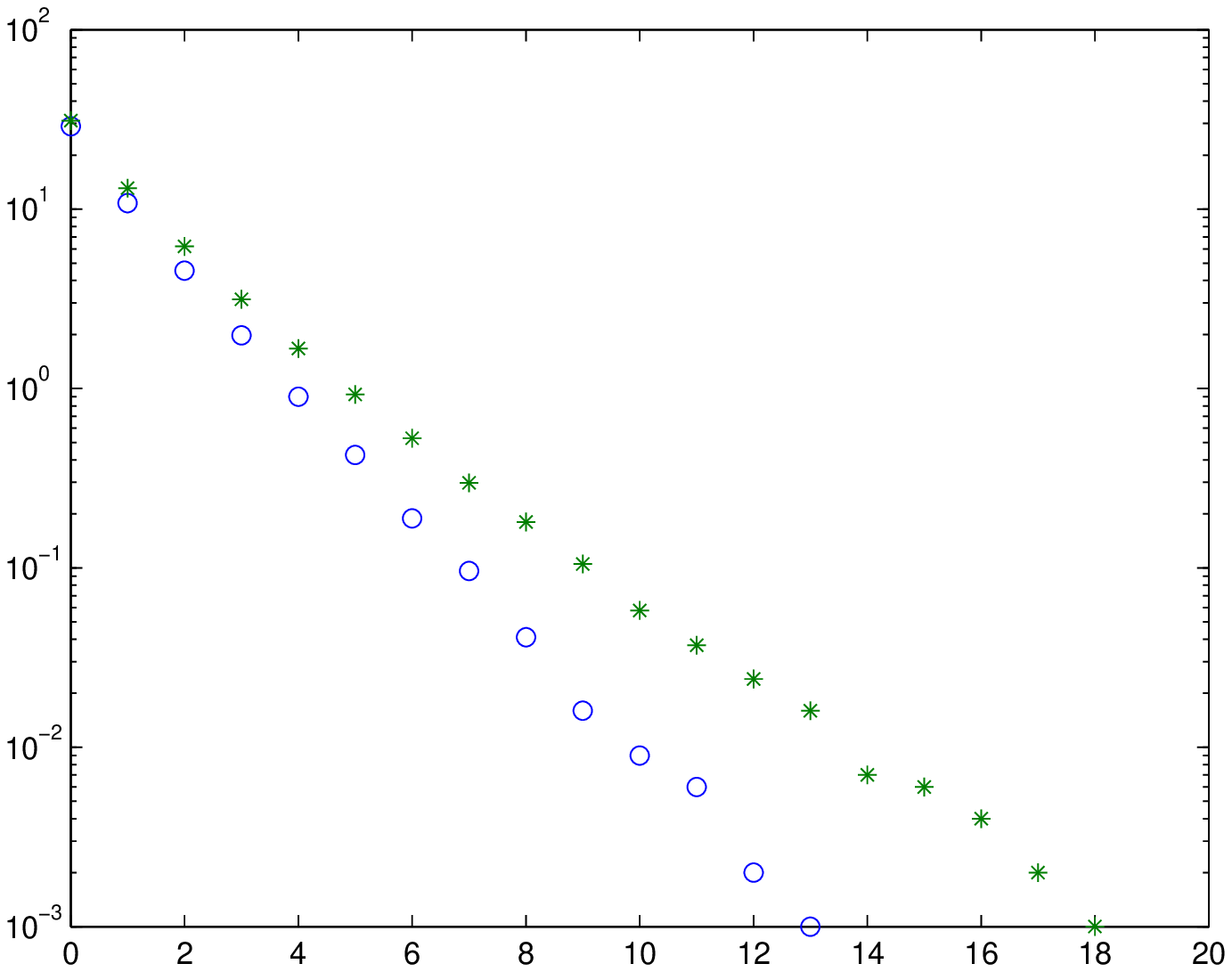,width=0.9\textwidth,height=0.55\textheight}}}
\end{center}
\caption{Simulation with $b(i)=i+1$ and $d(i)=0.5(i+1)$. The figure shows the tail probabilities of the fitness distribution (stars) and the in-degree distribution (circles).}
\end{figure}

\begin{figure}[p]
\begin{center}
\mbox{{\epsfig{file=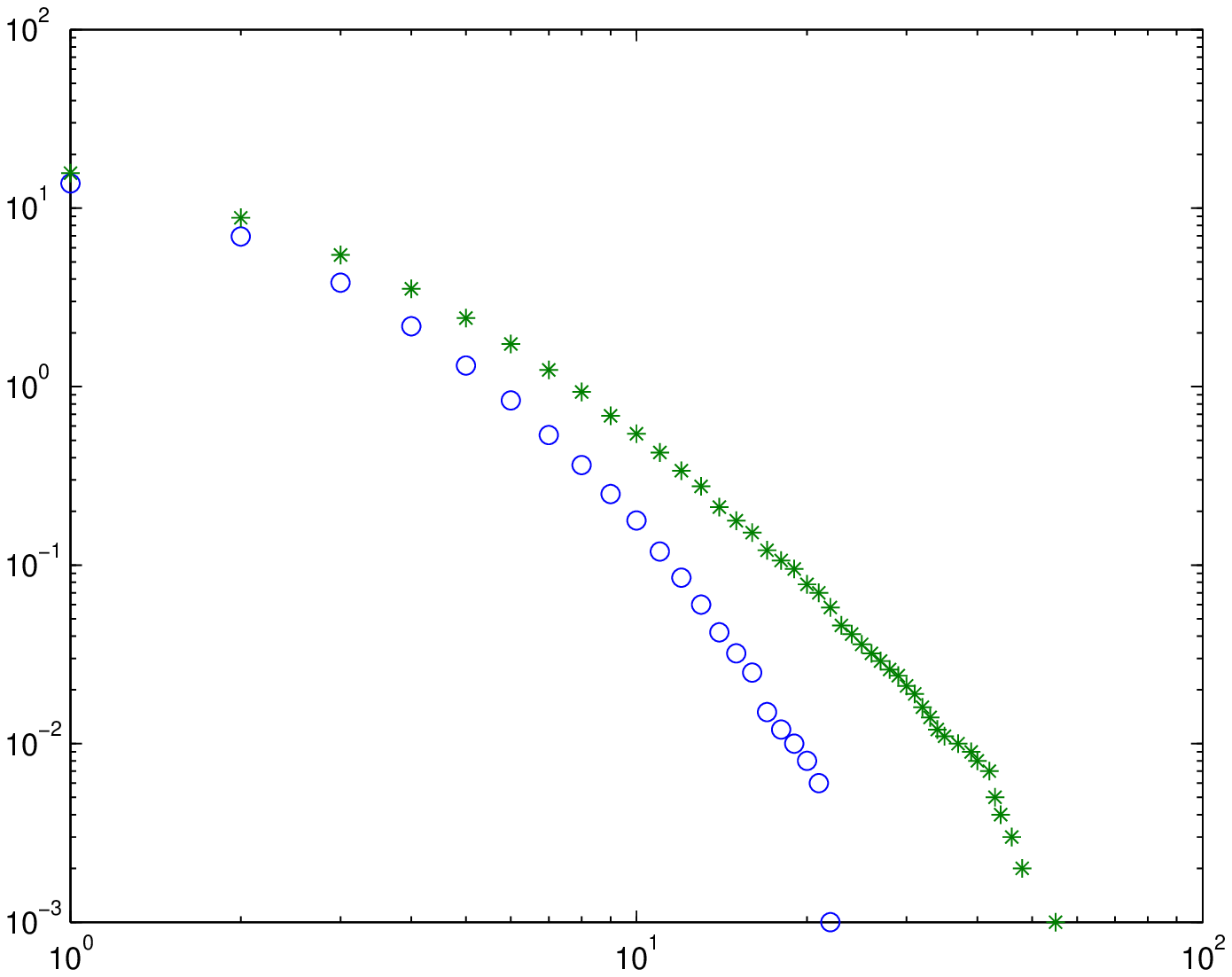,width=0.9\textwidth,height=0.55\textheight}}}
\end{center}
\caption{Simulation with $b(i)=i+1$ and $d(i)=0.5(i+1)$. The figure shows the tail probabilities of the fitness distribution (stars) and the in-degree distribution (circles).}
\end{figure}

\end{document}